\input amstex\documentstyle {amsppt}  
\pagewidth{12.5 cm}\pageheight{19 cm}\magnification\magstep1
\topmatter
\title Piecewise linear parametrization of canonical bases\endtitle
\author G. Lusztig\endauthor
\address Department of Mathematics, M.I.T., Cambridge, MA 02139\endaddress
\thanks Supported in part by the National Science Foundation\endthanks
\endtopmatter   
\document

\define\hcx{\hat{\cx}}

\define\uBB{\un{\BB}}

\define\us{\un s}

\define\ul{\un l}

\define\ur{\un r}

\define\uw{\un w}

\define\uI{\un I}

\define\uN{\un N}

\define\uW{\un W}
\define\ucx{\un{\cx}}

\define\ucb{\un{\cb}}
\define\utcx{\un{\tcx}}

\define\bi{\bar i}
\define\bj{\bar j}
\define\bs{\bar s}
\define\bn{\bar n}

\define\si{\sim}

\define\sqc{\sqcup}

\define\tcx{\ti\cx}

\define\lb{\linebreak}

\define\part{\partial}
\define\em{\emptyset}

\define\n{\notin}

\define\m{\mapsto}
\define\do{\dots}

\define\lra{\leftrightarrow}

\define\sub{\subset}    

\define\T{\times}
\define\ti{\tilde}
\define\nl{\newline}
\redefine\i{^{-1}}
\define\fra{\frac}
\define\un{\underline}
\define\ov{\overline}

\define\a{\alpha}
\redefine\b{\beta}

\define\g{\gamma}
\redefine\d{\delta}
\define\e{\epsilon}
\define\et{\eta}
\define\io{\iota}

\define\r{\rho}
\define\s{\sigma}

\define\th{\theta}

\redefine\l{\lambda}
\define\z{\zeta}
\define\x{\xi}

\define\Ph{\Phi}

\define\boc{\bold c}

\define\BB{\bold B}

\define\NN{\bold N}

\define\QQ{\bold Q}

\define\ZZ{\bold Z}

\define\cb{\Cal B}

\define\cx{\Cal X}

\define\fc{\frak c}
\define\fd{\frak d}

\define\BZ{BZ}
\define\KA{Ka}
\define\CB{L1}
\define\CBB{L2}
\define\COR{L3}
\define\QG{L4}
\define\PO{L5}
\head Introduction\endhead   
In \cite{\CB} the author introduced the canonical basis for the plus part of a quantized enveloping algebra of 
type $A,D$ or $E$. (The same method applies for nonsimplylaced types, see \cite{\COR, 12.1}.) Another approach to
the canonical basis was later found in \cite{\KA}. 
In \cite{\CB} we have also found that the set parametrizing the canonical basis has a natural piecewise linear
structure that is, a collection of bijections with $\NN^N$ such that any two of these bijections differ by
composition with a piecewise linear automorphism of $\NN^N$ (an automorphism which can be expressed purely in
terms of operations of the form $a+b,a-b,\min(a,b)$).
This led to the first purely combinatorial formula (involving only counting) for
the dimension of a weight space of an irreducible finite dimensional representation \cite{\CB} or the dimension of
the space of coinvariants in a triple tensor product \cite{\CBB, 6.5(f)}. (Later, different formulas in the same 
spirit were obtained by Littelman.) The construction of an analogous piecewise linear structure for the canonical
basis in the nonsimplylaced case (based on a reduction to the simplylaced case) was only sketched in \cite{\COR} 
partly because it involved an assertion whose proof only appeared later (in \cite{\QG, 14.4.9}): as Berenstein and
Zelevinsky write in \cite{\BZ, Proof of Theorem 5.2}, "Lusztig (implicitly)
claims that the transition map $R^{1212}_{2121}$ for $B_2$ is obtained from the transition map 
$R^{213213}_{132132}$ for type $A_3$...". In this paper we fill the gap in \cite{\COR} by making use of 
\cite{\QG, 14.4.9} which gives a relation between the canonical basis for a nonsimplylaced type and the canonical
basis for a simplylaced type with a given (admissible) automorphism. At the same time we slightly extend 
\cite{\QG, 14.4.9} by allowing type $A_{2n}$ with its non-admissible involution.

As an application we show that the canonical basis has a natural monoid structure and we define certain
"Frobenius" endomorphisms of this monoid.

\head 1. Parametrization\endhead
\subhead 1.1\endsubhead
In this paper a Cartan datum is understood to be a pair $(I,\cdot)$ where $I$ is a finite set and $(i\cdot j)$ is
a symmetric positive definite matrix of integers indexed by $I\T I$ such that 

$i\cdot i\in2\NN_{>0}$ for any $i\in I$;

$2\fra{i\cdot j}{i\cdot i}\in-\NN$ for any $i\ne j$ in $I$;
\nl
We say that $(I,\cdot)$ as above is 

-simply laced if $i\cdot j\in\{0,-1\}$ for any $i\ne j$ in $I$ and $i\cdot i=2$ for any $i\in I$; 

-irreducible if $I\ne\em$ and there is no partition $I=I'\sqc I''$ with  $I'\ne\em$, $I''\ne\em$, $i'\cdot i''=0$
for all $i'\in I'$, $i''\in I''$.

Let $(I,\cdot)$ be a Cartan datum. For $i\ne j$ in $I$ we have
$\fra{2i\cdot j}{i\cdot i}\fra{2j\cdot i}{j\cdot j}=0,1,2$ or $3$; accordingly, we set $h(i,j)=2,3,4$ or $6$. The
Weyl group $W$ of $(I,\cdot)$ is the group with generators $s_i(i\in I)$ and relations $s_i^2=1$ for $i\in I$ and
$s_is_js_i\do=s_is_js_i\do$ (both products have $h(i,j)$ factors) for $i\ne j$ in $I$. Let $l:W@>>>\NN$ be the
standard legth function relative to the generators $s_i$. Let $w_0$ be the unique element of $W$ such that
$l(w_0)$ is maximal. Let $N=l(w_0)$ and let $\cx$ be the set of sequences $i_*=(i_1,i_2,\do,i_N)$ in $I$ such that
$s_{i_1}s_{i_2}\do s_{i_N}=w_0$ (in $W$). We regard $\cx$ as the set of vertices of a graph in which $i_*,i'_*$ 
are joined if the sequences $i_*,i'_*$ coincide except at the places $k,k+1,k+2,\do,k+r-1$ where 

$(i_k,i_{k+1},\do,i_{k+r-1})=(p,p',p,\do)$, $(i'_k,i'_{k+1},\do,i'_{k+r-1})=(p',p,p',\do)$,
\nl
with $p\ne p'$ in $I$, $h(p,p')=r$. By a theorem of Iwahori and Tits,

(a) this graph is connected.

\subhead 1.2\endsubhead
Let $I=\{1,2,\do,2n\}$, $n\ge1$. For $i,j\in I$ we set $i\cdot j=2$ if $i=j$, $i\cdot j=-1$ if $i-j=\pm1$ and 
$i\cdot j=0$ otherwise. Then $(I,\cdot)$ is a simply laced irreducible Cartan datum. Define a permutation 
$\s:I@>>>I$ by $\s(i)=2n+1-i$ for all $i$. We have $\s(i)\cdot\s(j)=i\cdot j$ for any $i,j$ in $I$. 

\subhead 1.3\endsubhead
Let $I=\{1,2,\do,n-1,n,n'\}$, $n\ge1$. For $i,j\in[1,n]$ we set $i\cdot j=2$ if $i=j$, $i\cdot j=-1$ if $i-j=\pm1$
and $i\cdot j=0$ otherwise; we also set $n'\cdot n'=2$, $(n-1)\cdot n'=n'\cdot(n-1)=-1$, $i\cdot n'=n'\cdot i=0$ 
if $i<n-1$ or if $i=n$. Then $(I,\cdot)$ is a simply laced irreducible Cartan datum. It is irreducible if 
$n\ge2$. Define a permutation $\s:I@>>>I$ by $\s(i)=i$ for $i\in[1,n-1]$, $\s(n)=n'$, $\s(n')=n$. We have 
$\s(i)\cdot\s(j)=i\cdot j$ for any $i,j$ in $I$. 

\subhead 1.4\endsubhead
Let $\uI=\{\bar1,\bar2,\do,\bn\}$, $n\ge1$. For $i,j\in[1,n-1]$ we set $\bi\circ\bj=2$ if $i=j$, $\bi\circ\bj=-1$ if 
$i-j=\pm1$ and $\bi\circ\bj=0$ otherwise; we also set $\bn\circ\bn=4$, $\ov{n-1}\circ\bn=\bn\circ\ov{n-1}=-2$, 
$\bi\circ\bn=\bn\circ\bi=0$ if $i<n-1$. Then $(\uI,\circ)$ is an irreducible Cartan datum. 

\subhead 1.5\endsubhead
Let $(I,\cdot)$ be a simply laced Cartan datum and let $\s:I@>>>I$ be a permutation such that 
$\s(i)\cdot\s(j)=i\cdot j$ for any $i,j$ in $I$. Let $\uI$ be the set of orbits of $\s$ on $I$. For $\et\in\uI$ we
set $\d_\et=1$ if $\s(i)\cdot\s(j)=0$ for any $i\ne j$ in $\et$ and $\d_\et=2$, otherwise. We set 
$\d=\max_{\et\in\uI}\d_\et\in\{1,2\}$. We will assume that 

(a) either $\d=1$ or $(I,\cdot)$ is irreducible.
\nl
For any $\et\in\uI$ we set $\et\circ\et=2\d\i\d_\et|\et|$. For any $\et\ne\et'$ in $\uI$ we set 

$\et\circ\et'=-\d\i\d_\et\d_{\et'}|\{(i,j)\in\et\T\et';i\cdot j\ne0\}|$.
\nl
We show that $(\uI,\circ)$ is a Cartan 
datum. Assume first that $\d=1$. Let $\{x_\et;\et\in\uI\}$ be a collection of real numbers, not all zero. Let
$m=\sum_{\et,\et'\in\uI}x_\et x_{\et'}\et\circ\et'$. it is enough to show that $m>0$. For $i\in I$ let $y_i=x_\et$
where $i\in\et$. From the definitions we have $m=\sum_{i,i'\in I}y_iy_{i'}i\cdot i'$ and this is $>0$ since 
$(i\cdot i')$ is positive definite. Assume next that $\d=2$. We can assume that $(I,\cdot),\s$ are as in 1.2.
Denoting by $\bi$ the subset $\{i,2n+1-i\}$ of $I$ ($i\in[1,n]$) we see that $(\uI,\circ)$ is the same as that in
1.4 hence it is a Cartan datum.

\subhead 1.6\endsubhead
Let $(I,\cdot),\s$ be as in 1.3. We define $(\uI,\circ)$ in terms of $(I,\cdot),\s$ as in 1.5. Denoting by $\bi$ 
the subset $\{i\}$ of $I$ ($i\in[1,n-1]$) and the subset $\{n,n'\}$ if $i=n$, we see that $(\uI,\circ)$ is the 
same as that in 1.4.

\subhead 1.7\endsubhead
Let $(I,\cdot)$ be a simply laced Cartan datum. Let $W,l,w_0,N$ be as in 1.1. Let $K$ be either:

(i) a subgroup of the multiplicative group of a field which is closed under addition in that field;

(ii) a set with a given bijection $\io:\ZZ@>\si>>K$ with the operations $a+b,ab,a/b$ (on $K$) obtained by 
transporting to $K$ the operations $\min(a,b),a+b,a-b$ on $\ZZ$;

(iii) the subset $\io(\NN)$ of the set in (ii); this is stable under the operations $a+b,ab,a/(a+b)$.
\nl
Note that operations on $K$ in (i) and (ii) have similar properties; they are both examples of "semifields". (See
http://en.wikipedia.org/wiki/semifield.) The analogy between $K$ in (ii) and $K$ in (i) has been pointed out in 
\cite{\QG, 42.2.7} in connection with observing the analogy of the combinatorics of canonical bases and the 
geometry involved in total positivity.

Let $\tcx$ be the set of all objects $i_1^{c_1}i_2^{c_2}\do i_N^{c_N}$ (also denoted by $i_*^{c_*}$) where 
$i_*=(i_1,i_2,\do,i_N)\in\cx$, $c_*=(c_1,c_2,\do,c_N)\in K^N$. We regard $\tcx$ as the set of vertices of a graph
in which two vertices $i_*^{c_*}$, $i'_*{}^{c'_*}$ are joined if either

the sequences $i_*,i'_*$ coincide except at two places $k,k+1$ where $i'_k=i_{k+1},i'_{k+1}=i_k$ and
$i_k\cdot i_{k+1}=0$; the sequences $c_*,c'_*$ coincide except at the places $k,k+1$ where 
$c'_k=c_{k+1},c'_{k+1}=c_k$; or

the sequences $i_*,i'_*$ coincide except at three places $k,k+1,k+2$ where 

$(i_k,i_{k+1},i_{k+2})=(p,p',p)$, $(i'_k,i'_{k+1},i'_{k+2})=(p',p,p')$, 
\nl
with $p\cdot p'=-1$; the sequences $c_*,c'_*$ coincide except at the places $k,k+1,k+2$ where 

$(c_k,c_{k+1},c_{k+2})=(x,y,z)$, $(c'_k,c'_{k+1},c'_{k+2})=(x',y',z')$ 
\nl
with $x'=yz/(x+z)$, $y'=x+z$, $z'=xy/(x+z)$ or equivalently $x=y'z'/(x'+z')$, $y=x'+z'$, $z=x'y'/(x'+z')$.
\nl
We shall write $R_{i_*}^{i'_*}(c_*)=c'_*$ whenever $i_*^{c_*}$, $i'_*{}^{c'_*}$ are joined in the graph $\tcx$.
Then $R_{i_*}^{i'_*}:K^N@>>>K^N$ can be viewed as a bijection defined whenever $i_*,i'_*$ are joined in the graph
$\cx$.

Let $\cb$ be the set of connected components of the graph $\tcx$. For any $i_*\in\cx$ we define 
$\a_{i_*}:K^N@>>>\cb$ by $c_*\m\text{connected component of }i_*^{c_*}$. Note that:

(a) $\a_{i_*}$ is a bijection. 
\nl
If $K$ is as in 1.7(i) this follows from the proof of \cite{\QG, 42.2.4}. If $K$ is as in 1.7(ii) then, as in 
\cite{\QG, 42.2.7}, it can be viewed as a homomorphic image of a $K$ as in 1.7(ii) so that (a) holds in this case.
The case where $K$ is as in 1.7(iii) follows immediately from the case where $K$ is as in 1.7(ii), or it can be 
obtained directly from \cite{\QG, 42.1.9}.

For any $i_*,i'_*$ in $\cx$ we define a bijection $R_{i_*}^{i'_*}:K^N@>\si>>K^N$ as the composition
$R_{i^t_*}^{i^{t-1}_*}\do R_{i^1_*}^{i^2_*}R_{i^0_*}^{i^1_*}$ where $i^0_*,i^1_*,\do,i^t_*$ is a sequence of 
vertices in $\tcx$ such that $(i^0_*,i^1_*),(i^1_*,i^2_*),\do,(i^{t-1}_*,i^t_*)$ are edges of the graph $\cx$ and
$i^0_*=i_*$, $i^t_*=i'_*$ (such a sequence exists by 1.1(a)). This agrees with the earlier definition in the case
where $i_*,i'_*$ are joined in $\cx$. Note that $R_{i_*}^{i'_*}$ is independent of the choice above; it is equal 
to $\a_{i'_*}\i\a_{i_*}$.

\subhead 1.8\endsubhead
Let $(I,\cdot),\s$ be as in 1.5. Define $(\uI,\circ)$ as in 1.5. Define 
$W,l,w_0,N,\cx$ as in 1.1. Let $\uW,\ul,\uw_0,\uN,\ucx$ be the analogous objects defined in terms of 
$(\uI,\circ)$. The generators of $W$ are denoted by $s_i(i\in I)$ as in 1.1; similarly let $\us_\et(\et\in\uI)$ be
the generators of $\uW$. For any $\et\in\uI$ let $w_\et\in W$ be the longest element in the subgroup of $W$ 
generated by $\{s_i;i\in\et\}$; let $N_\et=l(w_\et)$. We can identify $\uW$ with the subgroup of $W$ 
generated by $\{w_\et;\et\in\uI\}$ by sending $\us_\et$ to $w_\et$. Then $\uw_0=w_0$ and $\ucx$ becomes the set of
sequences $\et_*=(\et_1,\et_2,\do,\et_{\uN})$ in $\uI^{\uN}$ such that $w_{\et_1}w_{\et_2}\do w_{\et_{\uN}}=w_0$.
We have $\uW=\{w\in W;\s(w)=w\}$ where $\s:W@>>>W$ is the automorphism given by $\s(s_i)=s_{\s(i)}$ for all $i$.
For any $\et\in\uI$ let $\cx^\et$ be the set of sequences $(h_1,h_2,\do,h_{N_\et})$ in $\et^{N_\et}$ such that 
$s_{h_1}s_{h_2}\do s_{h_{N_\et}}=w_\et$. 

Let $\tcx$ be as in 1.7. Let $\utcx$ be the set of all objects 
$\et_1^{\fc_1}\et_2^{\fc_2}\do\et_{\uN}^{\fc_{\uN}}$ (also denoted by $\et_*^{\fc_*}$) where 
$\et_*=(\et_1,\et_2,\do,\et_{\uN})\in\ucx$, $\fc_*=(\fc_1,\fc_2,\do,\fc_{\uN})\in K^{\uN}$. 

Let $\hcx$ be the set of all pairs $(\et_*^{\fc_*},\fd_*)$ where $\et_*^{\fc_*}\in\utcx$ and
$\fd_*=(\fd_1,\fd_2,\do,\fd_{\uN})$ is such that $\fd_j\in\cx^{\et_j}$ for $j\in[1,\uN]$. Let 
$(\et_*^{\fc_*},\fd_*)\in\hcx$. For $j\in[1,\uN]$, $k\in[1,N_j]$ (where $N_j=N_{\et_j}$) let $\e_{j,k}$ be the 
number of $k'\in[1,N_j]$ such that $h_{k'}=h_k$ where $\fd_j=(h_1,h_2,\do,h_{N_j})$. We have $\e_{j,k}\in\{1,2\}$.
Let $\e_j=\max_{k\in[1,N_j]}\d_{j,k}\in\{1,2\}$. Let 
$c^j_*=(\e_j\e_{j,1}\i\fc_j,\e_j\e_{j,2}\i\fc_j,\do,\e_j\e_{j,N_j}\i\fc_j)\in K^{N_j}$. Let 
$c_*=c^1_*c^2_*\do c^{\uN}_*\in K^N$ be the concatenation of $c^1_*,c^2_*,\do,c^{\uN}_*$. Let 
$i_*=\fd_1\fd_2\do\fd_{\uN}$ be the concatenation of $\fd_1,\fd_2,\do,\fd_{\uN}$. We have $i_*\in\cx$ and
$i_*^{c_*}\in\tcx$.

We show that the connected component of $i_*^{c_*}$ in $\tcx$ depends only on $\et_*^{\fc_*}$, not on $\fd_*$. Let 
$\fd'_*=(\fd'_1,\fd'_2,\do,\fd'_{\uN})$ be another sequence such that $\fd'_j\in\cx^{\et_j}$ for $j\in[1,\uN]$. 
Let $i'_*$ be the concatenation of $\fd'_1,\fd'_2,\do,\fd'_{\uN}$. Define $c'_*$ in terms of 
$(\et_*^{\fc_*},\fd'_*)$ in the same way as $c_*$ was defined in terms of $(\et_*^{\fc_*},\fd_*)$. We must show 
that $i_*^{c_*},i'_*{}^{c'_*}$ are in the same connected component of $\tcx$. We may assume that $I$ is a single
$\s$-orbit $\et$. Assume first that $\et=\{i,i'\}$ with $i\cdot i'=-1$. It is enough to show that 
$i^ci'{}^{2c}i^c$ and $i'{}^ci^{2c}i'{}^c$ are joined in $\tcx$ (where $c\in K$). This is clear since 
$c+c=2c,c(2c)/(c+c)=c$. Next assume that $\et$ is not of the form $\{i,i'\}$ with $i\cdot i'=-1$. Then 
$\et=\{i_1,i_2,\do,i_k\}$ where $s_{i_1},s_{i_2},\do,s_{i_k}$ commute with each other. It is enough to show that 
the connected component of $i_1^ci_2^c\do i_k^c$ in $\tcx$ does not depend on the order in which $i_1,i_2,\do,i_k$
are written (where $c\in K$); this is obvious.

We see that the map $\hcx@>>>\cb$ given by $(\et_*^{\fc_*},\fd_*)\m\text{connected component of }i_*^{c_*}$
factors through a map $s:\utcx@>>>\cb$ ($\cb$ as in 1.7).

We define a permutation $\s:\cx@>>>\cx$ by 

$i_*=(i_1,i_2,\do,i_N)\m(\s(i_1),\s(i_2),\do,\s(i_N))$
\nl
and a permutation $\s:\tcx@>>>\tcx$ by $i_*^{c_*}\m\s(i_*)^{c_*}$. This last permutation respects the graph
structure of $\tcx$ hence induces a permutation (denoted again by $\s$) of $\cb$.

We show that the image of $s:\utcx@>>>\cb$ is contained in the set $\cb^\s$ of fixed points of $\s:\cb@>>>\cb$.
Let $(\et_*^{\fc_*},\fd_*)\in\hcx$; we associate to it $i_*^{c_*}\in\tcx$ as above. For $j\in[1,\uN]$ we 
set $\fd'_j=(\s(h_1),\s(h_2),\do,\s(h_{N_j}))$ (where $\fd_j=(h_1,h_2,\do,h_{N_j})$, $N_j=N_{\et_j}$) and 
$\fd'_*=(\fd'_1,\fd'_2,\do,\fd'_{\uN})$. Let $i'_*$ be the concatenation of $\fd'_1,\fd'_2,\do,\fd'_{\uN}$. We 
have $i'_*\in\cx$. Now $i'_*{}^{c_*}$ is associated to $(\et_*^{\fc_*},\fd'_*)\in\hcx$ in the same way as 
$i_*{}^{c_*}$ is associated to $(\et_*^{\fc_*},\fd_*)\in\hcx$; hence $i_*^{c_*},i'_*{}^{c_*}$ are in the same 
connected component of $\tcx$ by an earlier argument. This verifies our claim.

Now let $\x\in\cb^\s$ and let $\et_*\in\ucx$. We show that $\x=s(\et_*^{\fc_*})$ for some $\fc_*\in K^{\uN}$. We 
can find $\fd_*=(\fd_1,\fd_2,\do,\fd_{\uN})$ such that $\fd_j\in\cx^{\et_j}$ for $j\in[1,\uN]$. Let 
$i_*=\fd_1\fd_2\do\fd_{\uN}$ be the concatenation of $\fd_1,\fd_2,\do,\fd_{\uN}$. We have $i_*\in\cx$ and by 
1.1(a) we can find $c_*\in K^N$ such that $i_*^{c_*}\in\x$. Let $\fd'_*$ be obtained from $\fd_*$ as in the 
previous paragraph and let $i'_*$ be the concatenation of $\fd'_1,\fd'_2,\do,\fd'_{\uN}$. We have 
$i'_*=\s(i_*)\in\cx$. Since $\x$ is $\s$-stable we see that $i_*^{c_*},i'_*{}^{c_*}$ are in the same connected 
component of $\tcx$. Now $c_*\in K^N$ can be viewed as the concatenation of $c^1_*,c^2_*,\do,c^{\uN}_*$ where
$c^j_*=(c^j_1,c^j_2,\do,c^j_{N_j})\in K^{N_j}$, $N_j=N_{\et_j}$. For $j\in[1,\uN]$ we write 
$\fd_j=(h_1,h_2,\do,h_{N_j})\in\et_j^{N_j}$ and we define 
$c'_*{}^j=(c'_1{}^j,c'_2{}^j,\do,c'_{N_j}{}^j)\in K^{N_j}$ by 

(i) $c'_k{}^j=c_{k'}^j$ where $\s(h_k)=h_{k'}$ if $s_{h_1},s_{h_2},\do,s_{h_{N_j}}$ commute with each other and

(ii) $c'_1{}^j=c_2^jc_3^j/(c_1^j+c_3^j),c'_2{}^j=c_1^j+c_3^j,c'_3{}^j=c_1^jc_2^j/(c_1^j+c_3^j)$ if 
$h_1\cdot h_2=-1$, $h_1=h_3$.   
\nl
Let $c'_*\in K^N$ be the concatenation of $c'_*{}^1,c'_*{}^2,\do,c'_*{}^{\uN}$. From the definitions we see that 
$i_*^{c_*},i'_*{}^{c'_*}$ are in the same connected component of $\tcx$. Hence $i'_*{}^{c_*},i'_*{}^{c'_*}$ are in
the same connected component of $\tcx$. Using the bijectivity of $\a_{i'_*}:K^N@>>>\cb$ (see 1.7(a)) we deduce 
that $c_*=c'_*$. Hence in (i) we have $c_k^j=c_{k'}^j$ whenever $\s(h_k)=h_{k'}$, hence $c_k^j$ is a constant 
$\fc_j$ when $k$ varies in $[1,N_j]$. Moreover in (ii) we have 
$c_1^j=c_2^jc_3^j/(c_1^j+c_3^j),c_2^j=c_1^j+c_3^j,c_3^j=c_1^jc_2^j/(c_1^j+c_3^j)$ hence 
$c_2^j=2\fc_j,c_1^j=c_3^j=\fc_j$ for some $\fc_j\in K$. Let $\fc_*=(\fc_1,\fc_2,\do,\fc_{\uN})\in K^{\uN}$. From 
the definitions we see that $s(\et_*^{\fc_*})$ is the connected component of $i_*^{c_*}$. Our claim is verified.

Assume that $\et_*\in\ucx$, $\fc_*\in K^{\uN}$, $\fc'_*\in K^{\uN}$ are such that $s(\et_*,\fc_*)=s(\et_*,\fc'_*)$. We 
show that $\fc_*=\fc'_*$. We can find $\fd_*=(\fd_1,\fd_2,\do,\fd_{\uN})$ such that $\fd_j\in\cx^{\et_j}$ for 
$j\in[1,\uN]$. We define $i_*^{c_*}\in\tcx$ in terms of $(\et_*^{\fc_*},\fd_*)$ as above and we define similarly
$i'_*{}^{c'_*}\in\tcx$ in terms of $(\et_*^{\fc'_*},\fd_*)$. Note that $i_*=i'_*$. By assumption,
$i_*^{c_*},i_*^{c'_*}$ are in the same connected component of $\tcx$. From 1.7(a) we see that $c_*=c'_*$. Now
$c_*\in K^N$ is the concatenation of $c^1_*,c^2_*,\do,c^{\uN}_*$ where
$c^j_*=(c^j_1,c^j_2,\do,c^j_{N_j})\in K^{N_j}$, $N_j=N_{\et_j}$. Similarly, $c'_*\in K^N$ is the concatenation of
$c'_*{}^1,c'_*{}^2,\do,c'_*{}^{\uN}$ where $c'_*{}^j=(c'_1{}^j,c'_2{}^j,\do,c'_{N_j}{}^j)\in K^{N_j}$ for
$j\in[1,\uN]$. We see that for any $j$ and any $k\in[1,N_j]$ we have $c^j_k=c'_k{}^j$. If $\e_j=1$ it follows that
$\fc_j=\fc'_j$. If $\e_j=2$ it follows that $(\fc_j,2\fc_j,\fc_j)=(\fc'_j,2\fc'_j,\fc'_j)$ hence again 
$\fc_j=\fc'_j$. We see that $\fc_*=\fc'_*$ as required.

From the previous two paragraphs we see that for any $\et_*\in\ucx$ the map $\a_{\et_*}:K^{\uN}@>>>\cb^\s$ given 
by $\fc_*\m s(\et_*^{\fc_*})$ is a bijection. 

For any $\et_*,\et'_*$ in $\ucx$ we define a bijection $R_{\et_*}^{\et'_*}:K^{\uN}@>>>K^{\uN}$ by
$R_{\et_*}^{\et'_*}=\a_{\et'_*}\i\a_{\et_*}$.

We regard $\utcx$ as the set of vertices of a graph in which two vertices $\et_*^{\fc_*}$, $\et'_*{}^{\fc'_*}$ are
joined if the sequences $\et_*,\et'_*$ coincide except at the places $k,k+1,k+2,\do,k+r-1$ where 

$(\et_k,\et_{k+1},\do,\et_{k+r-1})=(p,p',p,\do)$, $(\et'_k,\et'_{k+1},\do,\et'_{k+r-1})=(p',p,p',\do)$, 
\nl
with 
$p\ne p'$ in $\uI$, $h(p,p')=r$ and $R_{\et_*}^{\et'_*}(\fc_*)=\fc'_*$ or equivalently
$R_{\et'_*}^{\et_*}(\fc'_*)=\fc_*$. Here $h(p,p')$ is the analogue of $h(i,i')$ (see 1.1) for $(\uI,\circ)$ 
instead of $(I,\cdot)$.

Let $\ucb$ be the set of connected components of the graph $\utcx$. From the definitions we see that the map 
$s:\utcx@>>>\cb^\s$ factors through a map $\bs:\ucb@>>>\cb^\s$. We show that 

(a) $\bs$ is a bijection. 
\nl
The surjectivity of $\bs$ follows from the surjectivity of $s$. To prove that $\bs$ is injective we assume that 
$\et_*^{\fc_*}$, $\et'_*{}^{\fc'_*}$ are two elements of $\utcx$ such that 
$s(\et_*^{\fc_*})=s(\et'_*{}^{\fc'_*})$; we must show that $\et_*^{\fc_*}$, $\et'_*{}^{\fc'_*}$ are in the same 
connected component of $\utcx$. By the connectedness of the graph $\ucx$ (see 1.1(a)) we can find 
$\fc''_*\in K^{\uN}$ such that $\et'_*{}^{\fc'_*}$, $\et_*^{\fc''_*}$ are in the same connected component of 
$\utcx$. We have $s(\et'_*{}^{\fc'_*})=s(\et_*^{\fc''_*})$ hence $s(\et_*^{\fc''_*})=s(\et_*^{\fc_*})$. Using the
bijectivity of $\a_{\et_*}$ we deduce that $\fc_*=\fc''_*$. Thus, $\et'_*{}^{\fc'_*}$, $\et_*^{\fc_*}$ are in the
same connected component of $\utcx$ and our claim is verified.

Let $\et\in\uI$. We define a map $\un\l_\et:\ucb@>>>K$ by $\x\m\fc_1$ where $\et_*^{\fc_*}$ is any element of $\x$
such that $\et_1=\et$. (This map is well defined by an argument similar to that in \cite{\QG, 42.1.14}.) 
Similarly we define a map $\un\r_\et:\ucb@>>>K$ by $\x\m\fc_{\uN}$ where $\et_*^{\fc_*}$ is any element of $\x$ 
such that $\et_{\uN}=\et$. 

We define a map $\l_\et:\cb^\s@>>>K$ by $\x\m c_1$ where $i_*^{c_*}$ is any element of $\x$ such that $i_1\in\et$.
(This map is well defined.) Similarly we define a map $\r_\et:\cb^\s@>>>K$ by $\x\m c_N$ where $i_*^{c_*}$ is any
element of $\x$ such that $i_N\in\et$. 

From the definitions we have $\l_\et\bs=\un\l_\et$, $\r_\et\bs=\un\r_\et$, 

\subhead 1.9\endsubhead
We apply the definitions in 1.8 to $(I,\cdot),\s$ as in 1.3 and to $(\uI,\circ)$ as in 1.4, 1.6. Let 
$\et_*^{\fc_*},\et'_*{}^{\fc'_*}$ be two joined vertices of $\utcx$. We show:

(i) if $\et_*,\et'_*$ coincide except at the places $k,k+1$ where $(\et_k,\et_{k+1})=(\bi,\bi')$,
$(\et'_k,\et'_{k+1})=(\bi',\bi)$, $i-i'\n\{0,1,-1\}$ then $\fc_*,\fc'_*$ coincide except at the places $k,k+1$ 
where $(\fc_k,\fc_{k+1})=(x,y)$, $(\fc'_k,\fc'_{k+1})=(y,x)$;

(ii) if $\et_*,\et'_*$ coincide except at the places $k,k+1,k+2$ where
$(\et_k,\et_{k+1},\et_{k+2})=(\bi,\bi',\bi)$, $(\et'_k,\et'_{k+1},\et'_{k+2})=(\bi',\bi,\bi')$, ($i,i'$ in 
$[1,n-1]$, $i-i'=\pm1$, then $\fc_*,\fc'_*$ coincide except at the places $k,k+1,k+2$ where
$(\fc_k,\fc_{k+1},\fc_{k+2})=(x,y,z)$, $(\fc'_k,\fc'_{k+1},\fc'_{k+2})=(x',y',z')$
with $x'=yz/(x+z)$, $y'=x+z$, $z'=xy/(x+z)$ or equivalently $x=y'z'/(x'+z')$, $y=x'+z'$, $z=x'y'/(x'+z')$;

(iii) if $\et_*,\et'_*$ coincide except at the places $k,k+1,k+2,k+3$ where

$(\et_k,\et_{k+1},\et_{k+2},\et_{k+3})=(\ov{n-1},\bn,\ov{n-1},\bn)$,

$(\et'_k,\et'_{k+1},\et'_{k+2},\et'_{k+3})=(\bn,\ov{n-1},\bn,\ov{n-1})$
\nl
then $\fc_*,\fc'_*$ coincide except at the places $k,k+1,k+2,k+3$ where

$(\fc_k,\fc_{k+1},\fc_{k+2},\fc_{k+3})=(d,c,b,a)$, $(\fc'_k,\fc'_{k+1},\fc'_{k+2},\fc'_{k+3})=(d',c',b',a')$
\nl
and 

$d'=ab^2c/\e$, $c'=\e/\a$, $b'=\a^2/\e$, $a'=bcd/\a$ 

(or equivalently 

$d=a'b'{}^2c'/\e'$, $c=\e'/\a'$, $b=\a'{}^2/\e'$, $a=b'c'd'/\a'$) 
\nl
with the notation

$\a=ab+ad+cd$, $\e=ab^2+ad^2+cd^2+2abd$, 

$\a'=a'b'+a'd'+c'd'$, $\e'=a'b'{}^2+a'd'{}^2+c'd'{}^2+2a'b'd'$.
\nl
In case (i) and (ii) the result is obvious. In case (iii) we can assume that $n=2$ and we consider the sequence of
vertices of $\tcx$:
$$2^d2'{}^d1^c2'{}^b2^b1^a$$
$$2^d1^{\fra{bc}{b+d}}2'{}^{b+d}1^{\fra{cd}{b+d}}2^b1^a$$
$$2^d1^{\fra{bc}{b+d}}2'{}^{b+d}2^{\fra{ab(b+d)}{\a}}1^{\fra{\a}{b+d}}2^{\fra{bcd}{\a}}$$
$$2^d1^{\fra{bc}{b+d}}2^{\fra{ab(b+d)}{\a}}2'{}^{b+d}1^{\fra{\a}{b+d}}2^{\fra{bcd}{\a}}$$
$$1^{\fra{ab^2c}{\e}}2^{\fra{\e}{\a}}1^{\fra{dbc\a}{(b+d)\e}}2'{}^{b+d}1^{\fra{\a}{b+d}}2^{\fra{bcd}{\a}}$$
$$1^{\fra{ab^2c}{\e}}2^{\fra{\e}{\a}}2'{}^{\fra{\e}{\a}}1^{\fra{\a^2}{\e}}2'{}^{\fra{bcd}{\a}}2^{\fra{bcd}{\a}}$$
in which any two consecutive lines represent an edge in $\tcx$. This proves our claim.

Note that the expressions appearing in the coordinate transformation (iii) first appeared in the case 1.7(iii) in
a different but equivalent form in \cite{\COR, 12.5} and were later rewritten in the present form in 
\cite{\BZ, 7.1}. (In the last displayed formula in \cite{\COR, 12.5}, $a+d-f$ should be replaced by $c+d-f$.) In 
the cases 1.7(ii), 1.7(iii) the coordinate transformation $K^4@>>>K^4$ appearing in (iii) can be viewed as a coordinate 
transformation $\NN^4@>>>\NN^4$, $(d,c,b,a)\m(d',c',b',a')$, where

$d'=a+2b+c-\min(a+2b,a+2d,c+2d)$,

$c'=\min(a+2b,a+2d,c+2d)-\min(a+b,a+d,c+d)$,

$b'=2\min(a+b,a+d,c+d)-\min(a+2b,a+2d,c+2d)$,

$a'=b+c+d-\min(a+b,a+d,c+d)$,
\nl
since $a+b+d\ge\min(a+2b,a+2d)$.

\subhead 1.10\endsubhead
We apply the definitions in 1.8 to $(I,\cdot),\s$ as in 1.2. Then the associated $(\uI,\circ)$ is as in 1.4 (see 
1.5). Let $\et_*^{\fc_*},\et'_*{}^{\fc'_*}$ be two joined vertices of $\utcx$. We show that statements 
1.9(i)-(iii) hold in the present case. In case (i) and (ii) the result is obvious. In case (iii) we can assume 
that $n=2$ and we consider the sequence of vertices of $\tcx$:
$$1^a4^a2^b3^{2b}2^b1^c4^c2^d3^{2d}2^d$$
$$1^a2^b4^a3^{2b}2^b1^c4^c2^d3^{2d}2^d$$
$$1^a2^b4^a3^{2b}2^b1^c2^d4^c3^{2d}2^d$$
$$1^a2^b4^a3^{2b}1^{\fra{cd}{b+d}}2^{b+d}1^{\fra{bc}{b+d}}4^c3^{2d}2^d$$
$$1^a2^b4^a1^{\fra{cd}{b+d}}3^{2b}2^{b+d}1^{\fra{bc}{b+d}}4^c3^{2d}2^d$$
$$1^a2^b4^a1^{\fra{cd}{b+d}}3^{2b}2^{b+d}4^c1^{\fra{bc}{b+d}}3^{2d}2^d$$
$$1^a2^b4^a1^{\fra{cd}{b+d}}3^{2b}4^c2^{b+d}1^{\fra{bc}{b+d}}3^{2d}2^d$$
$$1^a2^b4^a1^{\fra{cd}{b+d}}3^{2b}4^c2^{b+d}3^{2d}1^{\fra{bc}{b+d}}2^d$$
$$1^a2^b1^{\fra{cd}{b+d}}4^a3^{2b}4^c2^{b+d}3^{2d}1^{\fra{bc}{b+d}}2^d$$
$$2^{\fra{bcd}{\a}}1^{\fra{\a}{b+d}}2^{\fra{ab(b+d)}{\a}}4^a3^{2b}4^c2^{b+d}3^{2d}1^{\fra{bc}{b+d}}2^d$$
$$2^{\fra{bcd}{\a}}1^{\fra{\a}{b+d}}2^{\fra{ab(b+d)}{\a}}3^{\fra{2bc}{a+c}}4^{a+c}3^{\fra{2ab}{a+c}}
T2^{b+d}3^{2d}1^{\fra{bc}{b+d}}2^d$$
$$2^{\fra{bcd}{\a}}1^{\fra{\a}{b+d}}2^{\fra{ab(b+d)}{\a}}3^{\fra{2bc}{a+c}}4^{a+c}
2^{\fra{d(b+d)(a+c)}{\a}}3^{\fra{2\a}{a+c}}2^{\fra{ab(b+d)}{\a}}1^{\fra{bc}{b+d}}2^d$$
$$2^{\fra{bcd}{\a}}1^{\fra{\a}{b+d}}2^{\fra{ab(b+d)}{\a}}3^{\fra{2bc}{a+c}}4^{a+c}
2^{\fra{d(b+d)(a+c)}{\a}}3^{\fra{2\a}{a+c}}1^{\fra{bcd\a}{\e}}2^{\fra{\e}{\a}}1^{\fra{ab^2c}{\e}}$$
$$2^{\fra{bcd}{\a}}1^{\fra{\a}{b+d}}2^{\fra{ab(b+d)}{\a}}3^{\fra{2bc}{a+c}}2^{\fra{d(b+d)(a+c)}{\a}}
4^{a+c}3^{\fra{2\a}{a+c}}1^{\fra{bcd\a}{\e}}2^{\fra{\e}{\a}}1^{\fra{ab^2c}{\e}}$$
$$2^{\fra{bcd}{\a}}1^{\fra{\a}{b+d}}3^{\fra{2bcd}{\a}}2^{b+d}3^{\fra{2ab^2c}{(a+c)\a}}4^{a+c}
3^{\fra{2\a}{a+c}}1^{\fra{bcd\a}{\e}}2^{\fra{\e}{\a}}1^{\fra{ab^2c}{\e}}$$
$$2^{\fra{bcd}{\a}}1^{\fra{\a}{b+d}}3^{\fra{2bcd}{\a}}2^{b+d}4^{\fra{\a^2}{\e}}3^{\fra{2\e}{\a}}
4^{\fra{ab^2c}{\e}}1^{\fra{bcd\a}{\e}}2^{\fra{\e}{\a}}1^{\fra{ab^2c}{\e}}$$
$$2^{\fra{bcd}{\a}}3^{\fra{2bcd}{\a}}1^{\fra{\a}{b+d}}2^{b+d}4^{\fra{\a^2}{\e}}3^{\fra{2\e}{\a}}
4^{\fra{ab^2c}{\e}}1^{\fra{bcd\a}{\e}}2^{\fra{\e}{\a}}1^{\fra{ab^2c}{\e}}$$
$$2^{\fra{bcd}{\a}}3^{\fra{2bcd}{\a}}1^{\fra{\a}{b+d}}2^{b+d}4^{\fra{\a^2}{\e}}3^{\fra{2\e}{\a}}
1^{\fra{bcd\a}{\e}}4^{\fra{ab^2c}{\e}}2^{\fra{\e}{\a}}1^{\fra{ab^2c}{\e}}$$
$$2^{\fra{bcd}{\a}}3^{\fra{2bcd}{\a}}1^{\fra{\a}{b+d}}2^{b+d}4^{\fra{\a^2}{\e}}1^{\fra{bcd\a}{\e}}
3^{\fra{2\e}{\a}}4^{\fra{ab^2c}{\e}}2^{\fra{\e}{\a}}1^{\fra{ab^2c}{\e}}$$
$$2^{\fra{bcd}{\a}}3^{\fra{2bcd}{\a}}1^{\fra{\a}{b+d}}2^{b+d}4^{\fra{\a^2}{\e}}1^{\fra{bcd\a}{\e}}
3^{\fra{2\e}{\a}}2^{\fra{\e}{\a}}4^{\fra{ab^2c}{\e}}1^{\fra{ab^2c}{\e}}$$
$$2^{\fra{bcd}{\a}}3^{\fra{2bcd}{\a}}1^{\fra{\a}{b+d}}2^{b+d}4^{\fra{\a^2}{\e}}1^{\fra{bcd\a}{\e}}
3^{\fra{2\e}{\a}}2^{\fra{\e}{\a}}1^{\fra{ab^2c}{\e}}4^{\fra{ab^2c}{\e}}$$
$$2^{\fra{bcd}{\a}}3^{\fra{2bcd}{\a}}1^{\fra{\a}{b+d}}2^{b+d}1^{\fra{bcd\a}{\e}}
4^{\fra{\a^2}{\e}}3^{\fra{2\e}{\a}}2^{\fra{\e}{\a}}1^{\fra{ab^2c}{\e}}4^{\fra{ab^2c}{\e}}$$
$$2^{\fra{bcd}{\a}}3^{\fra{2bcd}{\a}}2^{\fra{bcd}{\a}}1^{\fra{\a^2}{\e}}2^{\fra{\e}{\a}}
4^{\fra{\a^2}{\e}}3^{\fra{2\e}{\a}}2^{\fra{\e}{\a}}1^{\fra{ab^2c}{\e}}4^{\fra{ab^2c}{\e}}$$
$$2^{\fra{bcd}{\a}}3^{\fra{2bcd}{\a}}2^{\fra{bcd}{\a}}1^{\fra{\a^2}{\e}}4^{\fra{\a^2}{\e}}
2^{\fra{\e}{\a}}3^{\fra{2\e}{\a}}2^{\fra{\e}{\a}}1^{\fra{ab^2c}{\e}}4^{\fra{ab^2c}{\e}}$$
in which any two consecutive lines represent an edge in $\tcx$. This proves our claim.

\subhead 1.11\endsubhead
Define $\cb,\cb^\s$ as in 1.7, 1.8 in terms of $(I,\cdot),\s$ as in 1.2. The objects analogous to $(I,\cdot),\s,\cb,\cb^\s$ 
when $(I,\cdot),\s$ are taken as in 1.3 are denoted by \lb
$(I',\cdot),\s',\cb',\cb'{}^{\s'}$.

Let $\utcx$ be the graph attached to $(I,\cdot),\s$ as in 1.8 and let $\utcx'$ be the analogous graph attached to
$(I',\cdot),\s'$. From the results in 1.9, 1.10 we see that the graphs $\utcx,\utcx'$ are canonically isomorphic.
Hence the sets $\ucb$, $\ucb'$ of connected components of $\utcx,\utcx'$ are in canonical bijection. Combining 
with the canonical bijection $\ucb\lra\cb^\s$ (see 1.8(a)) and the analogous bijection $\ucb'\lra\cb'{}^{\s'}$ we
obtain a canonical bijection 

(a) $\cb^\s\lra\cb'{}^{\s'}$.

\subhead 1.12\endsubhead
In this subsection we take $K,\io$ as in 1.7(iii). Let $(I,\cdot)$ be a Cartan datum. Let $\frak f$ be the 
$\QQ$-algebra with $1$ with generators $\th_i (i\in I)$ and relations
$$\sum_{p,p'\in\NN;p+p'=1-2i\cdot j/(i\cdot i)}(-1)^{p'}(p!p'!)\i\th_i^p\th_j\th_i^{p'}=0$$
for $i\ne j$ in $I$. Let $\BB$ be the canonical basis of the $\QQ$-vector space $\frak f$ obtained by specializing
under $v=1$ the canonical basis of the quantum version of $\frak f$ defined in \cite{\CB,\QG}. For $i\in I$ and 
$b\in\BB$ we define $l_i(b)\in\NN$, by the requirement that $b\in\th_i^{l_i(b)}\frak f$, 
$b\n\th_i^{l_i(b)+1}\frak f$; we define $r_i(b)\in\NN$, by the requirement that $b\in\frak f\th_i^{l_i(b)}$, 
$b\n\frak f\th_i^{l_i(b)+1}$.

If $(I,\cdot)$ is simply laced and $\cb$ is as in 1.7 then we have a canonical bijection

(a) $\b:\BB@>\si>>\cb$
\nl
such that $\l_i\b=\io l_i$, $\r_i\b=\io r_i$ for all $i\in I$. (Here 
$\l_i,\r_i:\cb@>>>K$ are defined as $\l_\et,\r_\et$ in 1.8 in the case where
$\s=1$.) See \cite{\CB,\CBB}.

Now let $(I,\cdot),\s$ be as in 1.5. Let $(\uI,\circ)$ be as in 1.5. Let $\uBB$ be the analogue of $\BB$ when 
$(I,\cdot)$ is replaced by $(\uI,\circ)$ and let $\ul_\et:\uBB@>>>\NN$, $\ur_\et:\uBB@>>>\NN$ ($\et\in\uI$) be the
functions analogous to $l_i,r_i$ defined in terms of $(\uI,\circ)$. The algebra automorphism 
$\th_i\m\th_{\s(i)} (i\in I)$ of $\frak f$ restricts to a permutation of $\BB$ denoted again by $\s$. Let $\BB^\s$
be the fixed point set of $\s:\BB@>>>\BB$. For $\et\in\uI$ we define $l_\et:\BB^\s@>>>\NN$ by $l_\et(b)=l_i(b)$ 
with $i\in\et$; we define $r_\et:\BB^\s@>>>\NN$ by $r_\et(b)=r_i(b)$ with $i\in\et$.

We have the following result:

(b) {\it there is a canonical bijection $\g:\uBB@>\si>>\BB^\s$ such that $l_\et\g=\io\ul_\et$, $r_\et\g=\io\ur_\et$ for
any $\et\in\uI$.}
\nl
When $\d=1$ (see 1.5) this is established in \cite{\QG, 14.4.9}. Assume now that $\d=2$. Then $(I,\cdot),\s$ are 
as in 1.2. We shall use the notation in 1.11. Let $\BB',\s':\BB'@>>>\BB'$ be the analogues of $\BB,\s:\BB@>>>\BB$
when $(I,\cdot),\s$ are replaced by $(I',\cdot),\s'$. Since $(\uI,\circ)$ is the same when defined in terms of 
$(I,\cdot),\s$ or in terms of $(I',\cdot),\s'$ and since $\d=1$ for $(I',\cdot),\s'$ we see that we have a 
canonical bijection

(c) $\uBB\lra\BB'{}^{\s'}$.
\nl
We now consider the following composition of bijections

$\uBB\lra\BB'{}^{\s'}\lra\cb'{}^{\s'}\lra\cb^\s\lra\BB^\s$.
\nl
(The first bijection is given by (c). The fourth bijection is obtained from (a) which is compatible with the
actions of $\s$ by taking fixed point sets of $\s$. The second bijection is an analogue of the fourth bijection.
The third bijection is given by 1.11(a).) This bijection has the required properties. This establishes (b) in our
case.

\head 2. The "Frobenius" endomorphism $\Ph_e$ of $\uBB$\endhead
\subhead 2.1\endsubhead
We assume that we are in the setup of 1.8 and that $K,\io$ are as in 1.7(iii).
Following \cite{\PO, 9.11} we consider the monoid $M^+$ (with $1$) 
defined by the generators $\x_i^n$ ($i\in I,n\in\ZZ$) and the relations

(i) $\x_i^a\x_i^b=\x_i^{\min(a,b)}$ for any $i\in I$ and $a,b$ in $\ZZ$;

(ii) $\x_i^a\x_{i'}^b=\x_{i'}^b\x_i^a$ for any $i,i'\in I$ such that $i\cdot i'=0$ and any $a,b$ in $\ZZ$;

(iii) $\x_i^a\x_{i'}^b\x_i^c=\x_{i'}^{a'}\x_i^{b'}\x_{i'}^{c'}$ for any $i,i'$ in $I$ such that $i\cdot i'=-1$ and
any integers $a,b,c,a',b',c'$ such that $a'=b+c-\min(a,c)$, $b'=\min(a,c)$, $c'=a+b-\min(a,c)$, or equivalently
$a=b'+c'-\min(a',c')$, $b=\min(a',c')$, $c=a'+b'-\min(a',c')$.
\nl
(Here $\x_i^0$ is not assumed to be $1$.) 

{\it Remark.} In the last line of \cite{\PO, 9.9} one should replace "adding $c$ to the first entry of $\boc$" by
the text: "replacing the first entry $c_1$ of $\boc$ by $\min(c,c_1)$". In \cite{\PO, 9.10(a)}, $n+n'$ should be 
replaced by $\min(n,n')$. 

For any $i_*\in\cx$ we define a map $\z_{i_*}:K^N@>>>M^+$ by 

$c_*\m\x_{i_1}^{\io\i(c_1)}\x_{i_2}^{\io\i(c_2)}\do\x_{i_N}^{\io\i(c_N)}$.
\nl
From \cite{\PO, 9.10} we see that $\z_{i_*}$ is injective. Clearly its image is independent of the choice of 
$i_*$; we denote it by $M_0^+$. Note that $\x_i^nM_0^+\sub M_0^+$, $M_0^+\x_i^n\sub M_0^+$ for any $i\in I$,
$n\in\NN$. In particular, $M_0^+$ is a submonoid (without
$1$) of $M^+$. We define $\z:\tcx@>>>M_0^+$ by $i_*^{c_*}\m\z_{i_*}(c_*)$. This map is constant on any connected
connected component of $\tcx$ hence it induces a map $\bar\z:\cb@>>>M_0^+$ (necessarily a bijection). 

Now $\x_i^n\m\x_{\s(i)}^n$ (with $i\in I,n\in\ZZ$) defines a monoid automorphism $M^+@>>>M^+$ denoted again by
$\s$. It restricts to a monoid automorphism $M_0^+@>>>M_0^+$ denoted again by $\s$. This is compatible with the 
bijection $\s:\cb@>>>\cb$ via $\bar\z$. Note that the fixed points $M^{+\s},M^{+\s}_0$ are submonoids of
$M^+,M^+_0$. Consider the composite bijection $\uBB\lra\BB^\s\lra\cb^\s\lra M^{+\s}_0$. Here the first bijection 
is as in 1.12(b), the second bijection is induced by the one in 1.12(a) and the third bijection is induced by 
$\bar\z$. Via this bijection the monoid structure on $M_0^{+\s}$ becomes a monoid structure on $\uBB$.

\subhead 2.2\endsubhead
We show that the crystal graph structure on $\uBB$ introduced in \cite{\KA} is completely determined by the 
monoid structure of $M^{+\s}$. For simplicity we assume that $\s=1$ so that $\uBB=\BB$. We identify $\BB=M^+_0$
via $\bar\z$. As shown in \cite{\CBB} giving the crystal graph structure on $\BB$ is equivalent to giving for any
$i\in I$, $n\in\NN$ the subsets $l_i\i(n)$ (see 1.12) of $\BB$ and certain bijections $l_i\i(0)@>\si>>l_i\i(n)$.
Now $l_i\i(n)$ is exactly the set of $\x\in M_0^+$ such that $\x_i^a\x=\x$ for any $a\ge n$ and
$\x_i^a\x\ne\x$ for any $a\in[0,n-1]$. The inverse of the bijection $l_i\i(0)@>\si>>l_i\i(n)$ is given by 
$\x\m\x_i^0\x$.

\subhead 2.3\endsubhead
Let $\dot{\BB}$ the canonical basis of the "modified" enveloping algebra (see \cite{\QG, 31.1}) attached to a root
datum with Cartan datum $(I,\cdot)$. By combining \cite{\PO, 9.10} and the last line of \cite{\PO, 10.1} we obtain
a natural piecewise linear structure on $\dot{\BB}$ and a natural monoid structure on $\dot{\BB}$. 

\subhead 2.4\endsubhead
We now return to the setup in 2.1 and we fix an integer $e\ge1$. There is a well defined endomorphism 
$\Ph_e:M^+@>>>M^+$ (as a monoid with $1$) such that $\x_i^n\m\x_i^{en}$ for any $i\in I,n\in\ZZ$. 
This restricts to a monoid endomorphism $M^+_0@>>>M^+_0$. Moreover, it commutes with $\s:M^+@>>>M^+$ hence it 
restricts to a monoid endomorphism $M_0^{+\s}@>>>M_0^{+\s}$. Via the bijection $\uBB\lra M^{+\s}_0$ in 2.1 this 
becomes a monoid endomorphism $\uBB@>>>\uBB$ denoted again by $\Ph_e$. We call $\Ph_e$ the "Frobenius" 
endomorphism of the canonical basis $\uBB$. A similar definition applies to $\dot{\BB}$ instead of $\uBB$.


\widestnumber\key{BZ}
\Refs
\ref\key{\BZ}\by A.Berenstein and A.Zelevinsky\paper Tensor product multiplicities, canonical bases and totally 
positive varieties\jour Invent. Math.\vol143\yr2001\pages77-128\endref
\ref\key{\KA}\by M.Kashiwara\paper On crystal bases of the $q$-analogue of universal enveloping algebras\jour Duke
Math. J.\vol63\yr1991\pages465-516\endref
\ref\key{\CB}\by G.Lusztig\paper Canonical bases arising from quantized enveloping algebras\jour J. Amer. Math.
Soc.\vol3\yr1990\pages447-498\endref
\ref\key{\CBB}\by G.Lusztig\paper Canonical bases arising from quantized enveloping algebras, II, in "Common 
trends in mathematics and quantum field theories"\jour Progr. of Theor. Phys. Suppl.\vol102\yr1990\pages175-201
\endref
\ref\key{\COR}\by G.Lusztig\paper Introduction to quantized enveloping algebras\inbook New developments in Lie 
theory and their applications, ed. J.Tirao\bookinfo Progr. in Math.\vol105\publ Birkh\"auser Boston\yr1992\pages
49-65\endref
\ref\key{\QG}\by G.Lusztig\book Introduction to quantum groups\bookinfo Progr. in Math.\vol110\publ Birkh\"auser 
Boston\yr1993\endref
\ref\key{\PO}\by G.Lusztig\paper Total positivity in reductive groups\inbook Lie theory and geometry\bookinfo 
Progr.in Math.\vol123\publ Birkh\"auser Boston\yr1994\pages531-568\endref
\endRefs
\enddocument